\documentclass[times,doublespace]{anzsauth}

\usepackage{moreverb}
\usepackage{url}
\usepackage{grffile}
\usepackage{lineno}
\usepackage{lipsum}
\usepackage[UKenglish]{isodate}

\newcommand{\beq}{\begin{eqnarray*}}
\newcommand{\eeq}{\end{eqnarray*}}
\newcommand{\beqq}{\begin{equation}}
\newcommand{\eeqq}{\end{equation}}

\def\btheta{{\bvec \theta}}

\def\bSigma{{\bvec \Sigma}}
\def\bLambda{{\bvec \Lambda}}
\def\btheta{{\bvec \theta}}
\def\bmu{{\bvec \mu}}

\def\bvec#1{\mbox{\boldmath $#1$}}
\def\bU{{\bvec U}}

\def\bP{{\bvec P}}

\def\bZ{{\bvec Z}}
\def\bb{{\bvec b}}

\def\bM{{\bvec M}}
\def\bS{{\bvec S}}
\def\bJ{{\bvec J}}
\def\bT{{\bvec T}}
\def\bI{{\bvec I}}

\newcommand*\patchAmsMathEnvironmentForLineno[1]{%
  \expandafter\let\csname old#1\expandafter\endcsname\csname #1\endcsname
  \expandafter\let\csname oldend#1\expandafter\endcsname\csname end#1\endcsname
  \renewenvironment{#1}%
     {\linenomath\csname old#1\endcsname}%
     {\csname oldend#1\endcsname\endlinenomath}}%
     
\newcommand*\patchBothAmsMathEnvironmentsForLineno[1]{%
  \patchAmsMathEnvironmentForLineno{#1}%
  \patchAmsMathEnvironmentForLineno{#1*}}%
  
\AtBeginDocument{%
\patchBothAmsMathEnvironmentsForLineno{equation}%
\patchBothAmsMathEnvironmentsForLineno{align}%
\patchBothAmsMathEnvironmentsForLineno{flalign}%
\patchBothAmsMathEnvironmentsForLineno{alignat}%
\patchBothAmsMathEnvironmentsForLineno{gather}%
\patchBothAmsMathEnvironmentsForLineno{multline}%
}

\def\volumeyear{2018}

\begin{document}
\cleanlookdateon

\runningheads{How can the score test be consistent?}{N~KARAVARSAMIS, R~M~HUGGINS, AND B~J~T~MORGAN}
\title{How can the score test be consistent?}
\author{N. Karavarsamis\addressnum{1}\corrauth, G. Guillera-Arroita\addressnum{2}, R.M. Huggins\addressnum{1}, B.J.T. Morgan\addressnum{3}}
\affiliation{School of Mathematics and Statistics, The University of Melbourne}

\address{
\addressnum{1} School of Mathematics and Statistics, The University of Melbourne, Victoria 3010, Australia \\
\hspace*{1ex} Email: \texttt{nkarav$@$unimelb.edu.au} \\
\addressnum{2} School of Biosciences, The University of Melbourne, Victoria 3010, Australia \\
\addressnum{3} University of Kent, UK
}

\begin{abstract}
The score test statistic using the observed information is easy to compute numerically. Its large sample distribution under the null hypothesis is well known and is equivalent to that of the score test based on the expected information, the likelihood-ratio test and the Wald test. However, several authors have noted that under the alternative this no longer holds and in particular the statistic can take negative values. Here we examine the score test using the observed information in the context of comparing two binomial proportions under imperfect detection, a common problem in ecology when studying occurrence of species. We demonstrate through a combination of simulations and theoretical analysis that a new modified rule which we propose that rejects the null hypothesis when the observed score statistic is larger than the usual chi-square cut-off or is negative has power that is mostly greater to any other test. In addition consistency is largely restored. Our new test is easy to use and inference is always possible. Supplementary material for this article are available online as per journal instructions.
\end{abstract}


\ack{This work is funded by a Discovery Early Career Research Award from the Australian Research Council.
}

\maketitle
\section{Introduction}\label{sec-intro}

The score statistic computed using the observed information matrix has several practical advantages. It only requires the computation of estimates under the null hypothesis and the use of the observed rather
than the expected information does not require the computation of possibly complex expectations. A comparison of the observed and expected information matrices by \citet{efron78} for the one-parameter case
supported the use of the observed information over the Fisher (or expected) information in the likelihood-ratio test. They claim that the inverse of the observed information matrix leads to estimates of
variance which are closer to the data than those resulting from the Fisher, as well as agreeing more closely with Bayesian and even with fiducial analyses. They also make the point that the observed
information can be easier to compute, especially for complicated cases such as that given in \citet{cox58}. They give detailed evaluations of numerical results for moderate sample sizes, in addition to the
general asymptotic theory. In fact the observed score test statistic can be calculated using a numerical calculation of the observed information matrix, using difference approximation of differentials.
However, there are reservations on the use of the score statistic based on the observed information matrix as it has been observed that this version of the score statistic can be negative. Examples of
hypothesis tests that return negative score test statistics for the observed information matrix are given in \citet{lawrance87,yang03, godfrey01}. Situations that returned negative eigenvalues of the observed
information matrix are given in \citet{catchpole96,storer83,hosking84}. \citet{morgan07} examined the score statistic for the zero-inflated Poisson (ZIP) model and showed how it could be negative (Figure 1 of
\citet{morgan07}) when testing for zero inflation. Motivated by this result, \citet{freedman07} went on to further evaluate the use of the score statistic. \citet{verbeke07} summarise related issues with the score test.

We re-examine this problem and propose a solution in the context of comparing two binomial proportions under imperfect detection. This is a very common setting in ecology, when comparing the proportion of sites occupied by a species in
two regions or times \citep{guillera12}. Detection in wildlife surveys is generally imperfect, and species may remain completely undetected at sites they occupy. A common approach to account for detectability
is to conduct several survey visits to sites. The statistical model to describe these data is based on a zero-inflated Binomial (ZIB) distribution \citep{mac02}. The two-sample occupancy model involves the
probabilities of occupancy ($\psi_1$ and $\psi_2$) and detection ($p_1$ and $p_2$) of the species in each of the two regions (or times). The null hypothesis of interest is that the two occupancy probabilities
are the same: $H_0:\psi_1=\psi_2$. Thus under the null hypothesis there are three parameters and four under the alternative.

In this setting, we examine the observed score test statistic both with simulations and analytically.
We show that the commonly used rejection rule that rejects the null hypothesis if the value of the score statistic using the observed information matrix is greater than the usual $\chi^2_1$ cut-off can have power that is much lower than that of the likelihood-ratio test, the Wald test and the score test based on the expected information matrix.
However, in simulations we show that, if the rejection rule based on the observed score test statistic is modified to reject if its value is greater than the usual $\chi^2_1$ {\em or when it is negative}, the power is restored and mostly exceeds that of all the other tests. Our conclusion is that, with this new rejection rule, the score test based on the observed information matrix is indeed useful and can still be used, making inference always possible even when the statistic is negative. Under the null hypothesis the score statistic is positive so the new test does not alter the size of the test.

We organise our paper as follows.
In Section~\ref{sec-score}, we define the score test. In Section~\ref{sec-sims}, we conduct a simulation study to compare the likelihood-ratio test, the Wald test and the two versions of the score test. Our results confirm that the score test based on the observed information matrix can give negative test statistics and its naive use yields a test of low power, with power decreasing as the alternative moves away from the null. In Section \ref{sec-eigen} (and Appendix~\ref{sec-score2}) we examine the eigenvalues of the observed information matrix and its analytically computed expectation, under the alternative. We see that, as the alternative moves away from the null, the
 median over the simulations of the smallest eigenvalue of the observed information matrix becomes negative and so does the smallest eigenvalue of the expectation of
 this matrix.
This means that the observed information matrix is indefinite and hence quadratic forms can be positive or negative. This gives rise to negative values of the score statistic.
In Section~\ref{sec-modified}, we introduce our modified rejection rule and evaluate its performance in simulations.
 The paper ends with discussion in Section~\ref{sec-discuss}.

\section{The Score Test for the Species Occupancy Model}\label{sec-score}

In ecology, so-called `occupancy models' are used to estimate the probability $\psi$ that a species is present at a site (or equivalently, the proportion of the landscape that it occupies), while accounting
for imperfect detection \citep{mac02, guillera17}. To allow for detectability, detection/non-detection data are usually collected over repeated visits to the sites. A second, nuisance, parameter, the detection
probability $p$, is included in the model, and is the probability of detecting the species in a survey visit at a site where it is present.

In the standard model, detections are described as independent Bernoulli trials. Let $Y_i$ be the number of detections over $K$ visits at site $i$, $i=1,\dots,N$. Then
\begin{align}
Pr(Y_i = 0) &= 1 -  \psi + \psi(1-p)^K \nonumber \\
Pr( Y_i = y_i ) & = \psi p^{y_i}(1-p)^{K - y_i}, \quad y_i = 1,2,\ldots,K \quad i = 1,2,\ldots,N.
\end{align}
As the species is absent from some sites, the number of detections follows a zero-inflated binomial distribution (ZIB), with the level of zero-inflation set by $1-\psi$.
Let us now suppose that we want to compare species occupancy ($\psi$) in two samples. Without loss of generality, hereafter we assume these represent two geographical regions, but the same applies for comparison in time. We consider the two regions are labelled 1 and 2 with associated occupancy probabilities $\psi_1$, $\psi_2$ and detection probabilities $p_1$ and $p_2$ at sites within the regions.
First consider the likelihood for a single region $j$. The probability of detecting the species at least once at an occupied site from region $j$ is $\theta_j = 1 - (1-p_j)^{K_j}$. Let $s_{d_j}$ denote the number of sites where individuals were detected and $d_j$ the total number of detections, for region $j$. Note that for the single region case $d = \sum_{i=1}^N Y_i$.
Then, following \citet{guillera10}, the likelihood component for region $j$ is
\[ L_j = \{ \psi_j^{s_{d_j}} p_j^{d_j} (1-p_j)^{K_j s_{d_j} - d_j} \} (1- \psi_j \theta_j)^{N_j - s_{d_j}}, \quad j=1,2.\]
The score equations for region $j$ are given by
\[S_{j1}=s_{d_j}/ \psi_j - (N_j - s_{d_j})\theta_j/(1 - \psi_j \theta_j) = 0\] and
\[S_{j2}=d_j/p_j - (s_{d_j} K_j - d_j)/(1-p_j) - (N_j - s_{d_j})\psi_j K_j(1 - p_j)^{K_j-1}/(1 - \psi_j \theta_j) = 0, \] 
$j=1,2.$ \\
Let $\btheta=(\psi_1,p_1,\psi_2,p_2)^T$ be the vector of model parameters for the two-sample model.
 The likelihood for the two-sample model is the product of the two single-sample likelihoods assuming independence, and the resulting unconstrained score function is $\bS(\btheta)=(S_{11},S_{12},S_{21},S_{22})^T$.
 Also let $\bJ(\btheta)=\partial \bS(\btheta)/\partial\btheta^T=\bS'(\btheta)$ be the observed information matrix.  Standard likelihood theory gives the large-sample null distribution (i.e.\ $\chi^2_1$) of the score statistic using both the expected and observed information matrices.

Our interest is in the behaviour of the score test based on the observed information matrix, so we define this formally. The definitions of the likelihood-ratio, Wald test and
 the score test based on the expected information matrix are well known and omitted.
 Let
\[
     \bM=\begin{pmatrix}1&0&0\\0&1&0\\1&0&0\\0&0&1
     \end{pmatrix}.
\]

We denote the true parameter values by $\btheta_T= (\psi_{1T}, p_{1T}, \psi_{2T}, p_{2T})^T$. Consider $H_0:\psi_1=\psi_2=\psi$. Let $\btheta'=(\psi,p_1,p_2)^T$ and let $\bS_0(\btheta')=\bM^T\bS(\bM \btheta')$ be the score function under $H_0$. We suppose throughout that $\btheta'_S$ satisfies $E_{\btheta_T}(\bS_0(\btheta'_S))=0$ and the maximum likelihood estimate $\widehat\btheta'_S$ is a solution of $\bS_0(\btheta')=0$.
    Then the score statistic defined in terms of the observed information is
    \begin{equation}
    \bT_O(\widehat\btheta'_S)=\bS(\bM\widehat\btheta'_S)^T \bJ(\bM\widehat\btheta'_S)^{-1}\bS(\bM\widehat\btheta'_S), \label{eq-score-2}
    \end{equation}
and replacing $\bJ(\bM\widehat\btheta'_S)$ with $E_{\btheta_T}(\bM\bJ(\btheta'_S))$ evaluated at $\btheta_S'=\widehat \btheta_S'$ gives the expected score test statistic $\bT_E(\widehat\btheta'_S)$.
 In our setting, this will have a chi-square distribution with one degree of freedom under $H_0$ asymptotically.

\section{Performance of competing tests}\label{sec-sims}

In this section, we conduct a simulation study to compare the powers of the Wald test (probability scale), likelihood-ratio test and the two versions of the score test (based on the expected and observed information matrices).
We take $p_1=p_2=0.5$, $K_1=K_2=3$ visits and $N_1=N_2=50$ sites and refer throughout to this combination of parameter values as the `standard configuration'.
We consider $\psi_1=0.8$ and $\psi_2=(1-R)\psi_1$, where $R$ is the proportional decline in species occurrence. We sweep through values of $R \in 0,\dots,0.9$ in steps of 0.025. We run 50000 simulations at each value of $R$. Our simulation code is available as Supplementary Material. For illustrative purposes, we also show results for scenarios with $\psi_1=0.4$ in Appendix \ref{sup}. Results are displayed for simulations where the numerical optimization did not fail (for example when the optimization process did not converge to any value or the hessian was not invertible) and estimates produced were within the parameter space $(0,1)$ for $\psi$ and $p$.

\begin{figure}[h!]
\begin{center}
\includegraphics[width=.8\textwidth]{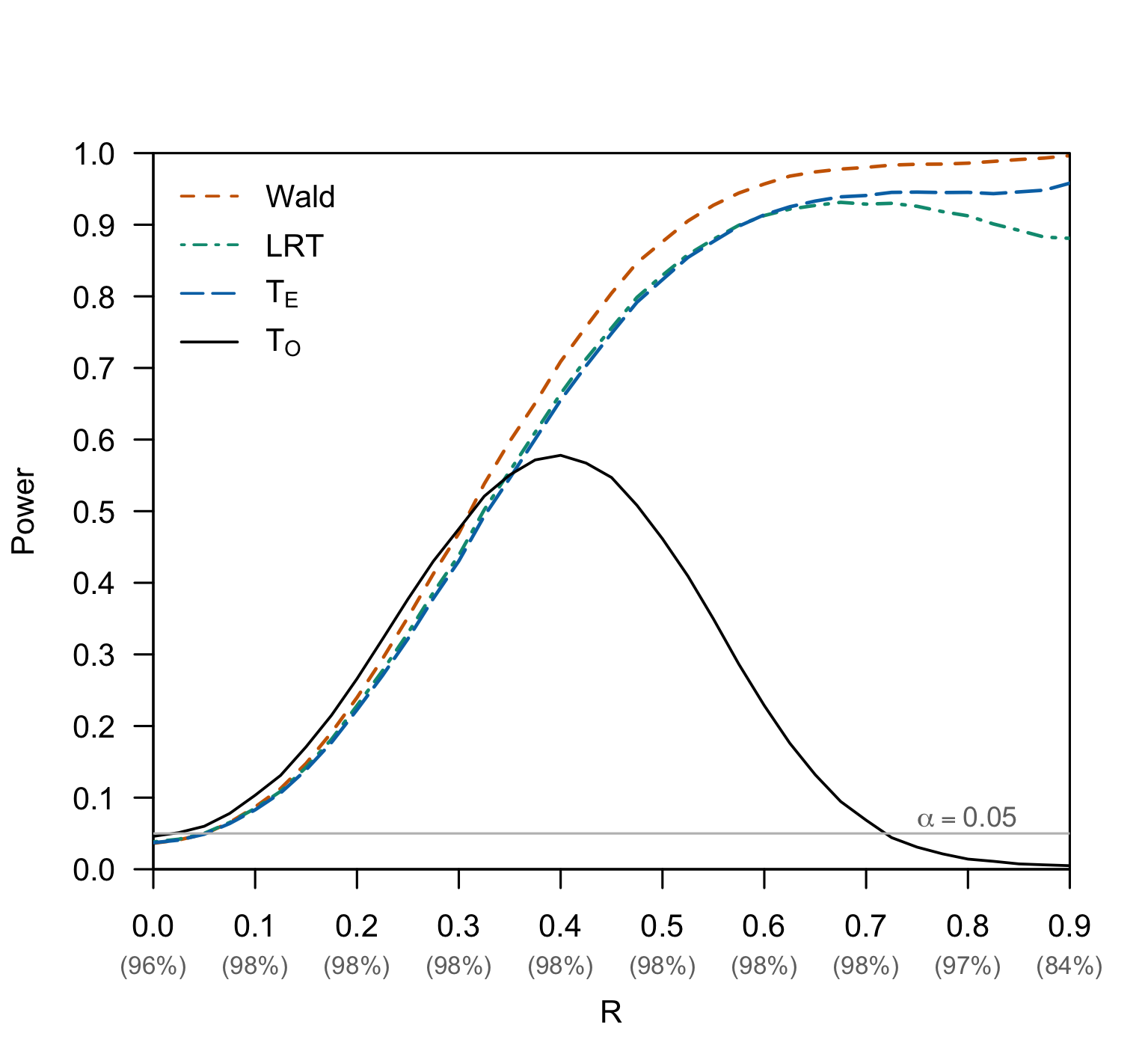}
\end{center}
\caption{Power of the four tests we consider (LRT = Likelihood-ratio test; Wald = Wald test; $T_E$ = Score test using expected information; $T_O$ = Score test using the observed information), for different effect sizes $R$ and computed via simulation (50000 datasets per $R$ point). Results are displayed for simulations where numerical optimization did not fail (percentages for $T_O$ are listed in parenthesis under the $x$-axis).}
\label{fig-power1}
\end{figure}

\begin{figure}[h!]
\begin{center}
\includegraphics[width=.8\textwidth]{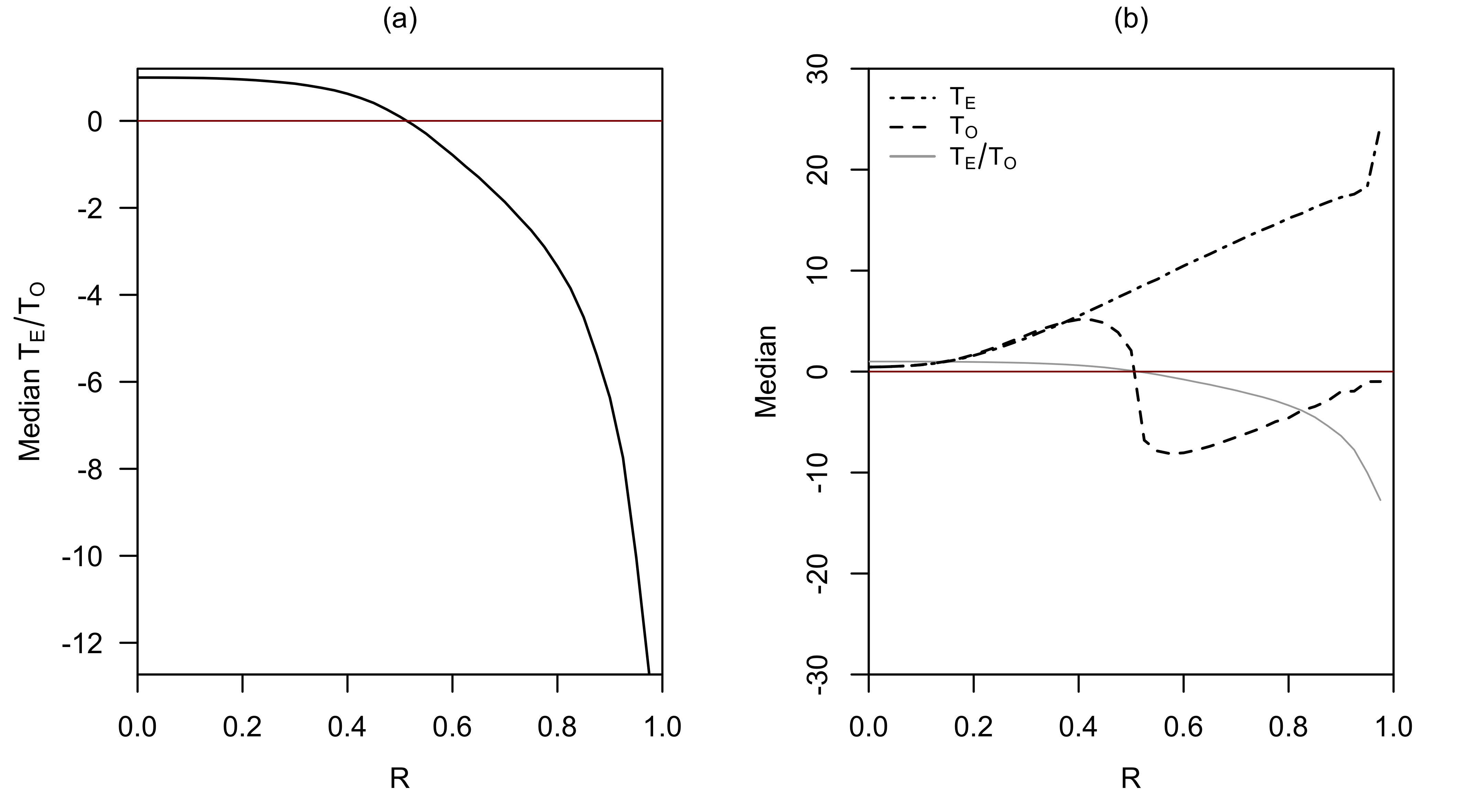}
\end{center}
\caption{Score test statistics versus effect size $R$ for data simulated under the standard configuration with $\psi_1=0.8$. In (a), median of the ratio between the expected score and the observed score test statistics ($T_E/T_O$). In (b), same ratio, together with the median of the expected score ($T_E$) and observed score ($T_O$) test statistics, plotted separately.}
\label{fig-2}
\end{figure}

The results (Figure~\ref{fig-power1}) show that the Wald test, the likelihood-ratio test (LRT) and the score test based on the expected information matrix ($T_E$) all have similar powers that increase as the effect size $R$ increases. However the score test from the observed information matrix ($T_O$) shows a dramatic decline in power after $R\approx 0.4$. Along the lines of Figure 1 of \citet{morgan07}, if we evaluate the relationship between the expected score test statistic and the observed, the median of their ratio becomes negative at 0.5 as $R$ increases (Figure~\ref{fig-2} (a)). At $\psi_1=\psi_2$, the null hypothesis is true with effect size equal to zero, i.e.\ $R=0$. Then the score statistics are equal and their ratio is exactly equal to 1, as seen in Figure~\ref{fig-2} (a).
When the medians of the simulated values of the score statistics are plotted separately (Figure~\ref{fig-2} (b)), we see that the median of the observed score test statistic reaches a maximum at 0.4, then declines
and later increases. It becomes negative at 0.5.  That is, at  $R\approx 0.5$ half of the values of the observed score statistic are positive and half are negative.
This explains why the ratio of the expected to observed is a declining function and that it becomes negative at 0.5 (Figure~\ref{fig-2} (a)). At $\approx 0.6$ the observed score statistic reaches another turning point and begins to increase however its median remains negative. 
Similar results and a decline at $\approx 0.5$ are obtained when we take $\psi_1=0.4$ (see Figure \ref{fig-power04}, Appendix \ref{sup}).

When we consider only those simulations where the observed score statistic is positive $(T_O^+)$, we find there is good agreement between the expected $(T_E)$ and observed $(T_O^+)$ score test, i.e. both accept or reject the null hypothesis for a given dataset (Table~\ref{tab-agree}(a)). The differences are predominantly where the test based on the observed score statistic rejects the null hypothesis and that based on the expected score statistic does not (see Figure \ref{fig-7}, Appendix \ref{sup}). The number of tests $n$ for $T_O^+$ is significantly reduced. As $R$ increases, the number of datasets with positive tests $n$ decreases substantially, e.g. when $R=0.8$ there are 1569 $(=n)$ positive tests of the 50000 simulated datasets (Table~\ref{tab-agree}(a)). We wish to increase $n$.
The observed score statistic has a size that exceeds slightly the significance level, however this is not due to negative values of the statistic.

\begin{table}[ht]
    \caption{Agreement (Agr.) for accepting/rejecting $H_0$ between tests based on the observed and expected score statistics, for data simulated under the standard configuration with $\psi_1=0.8$. In a) only datasets where the observed score statistic is positive ($T_O^+$) are considered. In b) both positive and negative score statistics are considered under our new modified rejection rule ($T_O^*$). In each case, $n$ is the number of simulated tests considered (out of the 50000 simulations). As effect size $R$ increases, there are fewer simulations with a positive observed score test statistic, hence the corresponding decrease in $n$.}
\vspace{2mm}
\centering
\begin{small}
\begin{tabular}{|l|l|rrrrrrrrrr|}
  \hline
 &$R$ & 0 & 0.1 & 0.2 & 0.3 & 0.4 & 0.5 & 0.6 & 0.7 & 0.8 & 0.9 \\
  \hline
a) $T_O^+$  & Agr.\ & 0.99 & 0.98 & 0.95 & 0.91 & 0.88 & 0.85 & 0.85 & 0.86 & 0.94 & 0.99 \\ \hline
  Positive &$n$ & 49974 & 49948 & 49574 & 47542 & 40502 & 26986 & 12388 & 3805 & 1569 & 5317 \\ \hline \hline
b) $T_O^*$  & Agr.\ & 0.99 & 0.98 & 0.95 & 0.91 & 0.90 & 0.91 & 0.94 & 0.95 & 0.95 & 0.96 \\
 Modified &$n$  & 49981 & 49992 & 50000 & 49999 & 49998 & 49992 & 49982 & 49982 & 49929 & 48454 \\
   \hline
\end{tabular}
\end{small}
\label{tab-agree}
\end{table}

\section{Eigenvalues}\label{sec-eigen}

\begin{figure}[ht]
\begin{center}
\includegraphics[width=.8\textwidth]{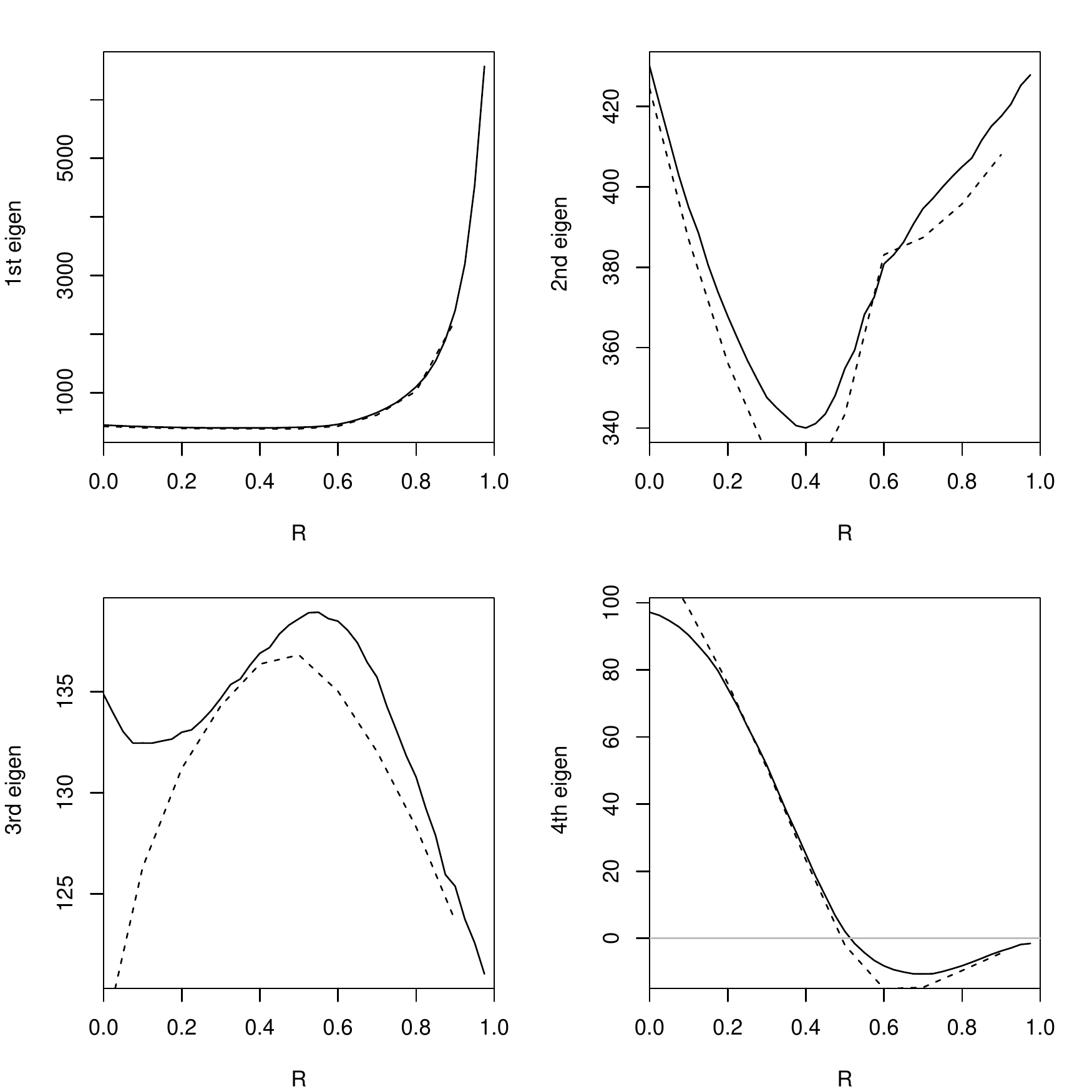}
\end{center}
\caption{Eigenvalues for the observed information matrix, for different values of effect size $R$ under the standard configuration with $\psi_1=0.8$. Solid lines are medians obtained from  simulations (50000 at each value of $R$). Dashed lines are eigenvalues of $E_{\btheta_T}(\bJ(\bM\btheta'_S))$.}
\label{fig-3}
\end{figure}

When we inspect the eigenvalues of the observed information matrix in our simulations, we find that, as $R$ increases, the median of the smallest eigenvalue becomes negative at $R \approx 0.5$ (Figure~\ref{fig-3}). To further investigate this result, we compute the expectation of the observed information matrix.
Now, as $\btheta_T$ is the true value, $\widehat\btheta' \buildrel P \over \longrightarrow \btheta'_S$ and  $E_{\btheta_T}(\bJ(\bM\btheta'_S))$ may be readily computed.  This requires computing $\btheta'_S$ for a given $\btheta_T$. Then for $j=1,2$, $E(s_{d_j}) = N_j \psi_{jT} \theta_{jT}$, $E(d_j \vert s_{d_j}) = K_j s_{d_j} p_{jT} / \theta_{jT}$ and $E(d_j) = K_j N_j \psi_{jT} p_{jT}$.
In Figure \ref{fig-3},
	we plot the eigenvalues of $E_{\btheta_T}(\bJ(\bM\btheta'_S))$ for the case $\psi_1=0.8$ under our standard configuration. We find that the smallest eigenvalue has a value of zero at $R \approx 0.5$.
  Despite the eigenvalues of $E_{\btheta_T}(\bJ(\bM\btheta'_S))$ being positive for $R<0.5$,  the estimator and the observed information matrix are random. Then for $R<0.5$  there will be some outcomes with negative eigenvalues leading to  negative values of the observed score statistic, as is apparent from Figure \ref{fig-2}. That is, the observed information matrix is indefinite in this case.

The turning point at $R = 0.4$ previously found in Figure~\ref{fig-2}(b) reflects the turning point at 0.4 for the 2nd eigenvalue, as seen in Figure~\ref{fig-3} for simulated and analytical results. In Appendix \ref{sec-score2} we examine the score statistic in more detail through the usual large-sample approximations.

\section{A Modified Rejection Rule}\label{sec-modified}

\begin{figure}[h!]
\begin{center}
\includegraphics[width=.8\textwidth]{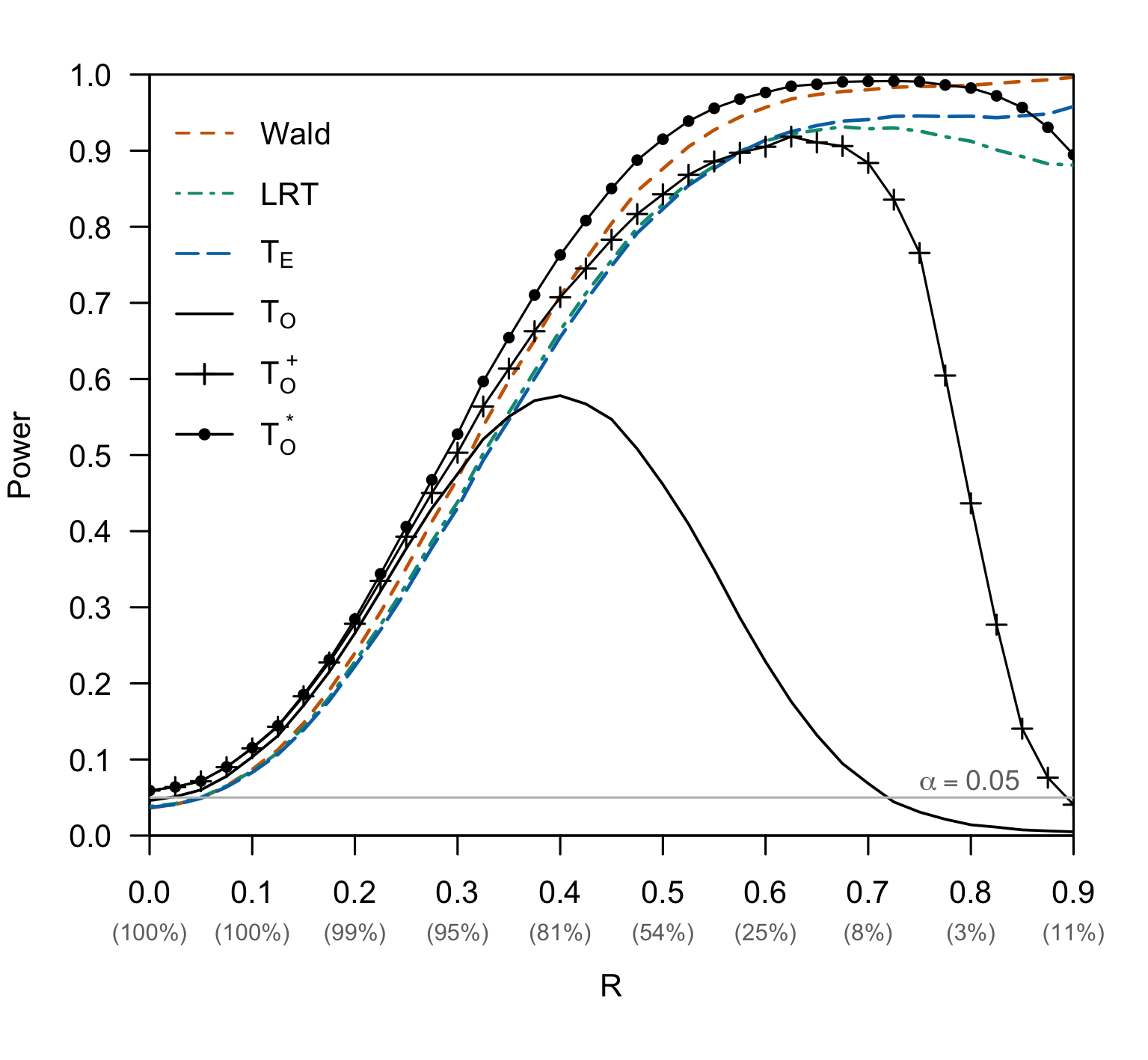}
\end{center}
\caption{Power plot for the new modified observed score test $T_O^*$, and for the observed score test when only the positive test statistics are considered $T_O^+$ (in addition to other tests as per Figure~\ref{fig-power1}). Results are displayed for simulations, from a total of 50000, where numerical optimization did not fail. Percentage of simulations that returned a positive observed score test ($T_O^+$) are listed in parenthesis under the $x$-axis.}
\label{fig-4}
\end{figure}

\begin{figure}[h!]
\begin{center}
\includegraphics[width=.8\textwidth]{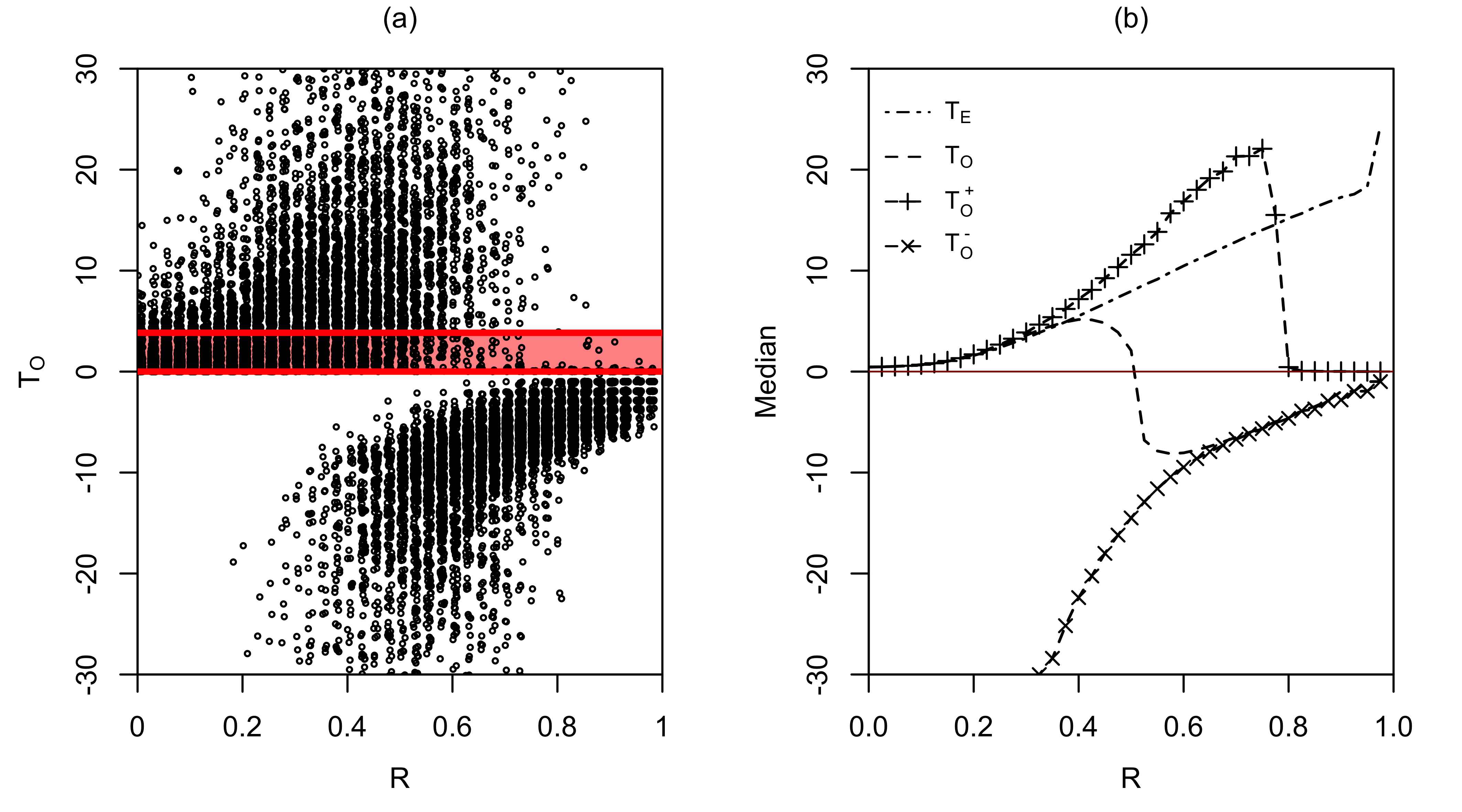}
\end{center}
\caption{Visual display of the new modified rejection rule for $\psi_1=0.8$ under the standard configuration. In (a) observed score test statistic $T_O$ versus effect size $R$ (for clarity, only 500 simulations are shown for each $R$). Horizontal lines at 0 and $\chi^2_1$ for $\alpha = 0.05$, i.e. at 3.814, give bounds for the acceptance region. Power for each $R$ is the proportion of simulations that lie outside the acceptance region. In (b), medians of positive ($T_O^+$), negative ($T_O^-$) and new modified rule ($T_O^*$) values of the observed score statistic, together with medians for all values ($T_O$) and for the expected score ($T_E$), obtained from 50000 simulations.}
\label{fig-8}
\end{figure}

The usual rejection rule  is to reject $H_0$ if the score statistic obtained using the observed information is larger than the chosen $\chi^2$ critical value. This is based on the null distribution where the observed information matrix is positive definite so that the observed score statistic is positive. However, we have seen that an indefinite information matrix can also be evidence that the null hypothesis is false (as demonstrated by the appearance of the negative 4th eigenvalue in Figure \ref{fig-3} as well from the negative values in Figures \ref{fig-5} and \ref{fig-6} in Appendix \ref{sec-score2}).
Because the observed information matrix may be indefinite under the alternative, its inverse gives some positive and negative eigenvalues. We therefore propose a new modified rejection rule based on rejecting according to being larger than  the usual $\chi^2$ critical value or  when the observed score test statistic is negative.
Our simulation results (Figure~\ref{fig-4}) clearly show that this new rejection rule, $T^*_O$, greatly improves the power of the observed score test. It mostly exceeds the power of the other tests, including the Wald and expected score test, over effect size $R$. Figure~\ref{fig-8} provides a visual display of the modified rejection rule.

Checking the percentage of simulations where the new rejection rule agrees with that of the score test based on the expected information matrix, we find that there is generally good agreement between the two approaches and $n$ is substantially increased for greater $R$ ($>$ 0.4) (Table \ref{tab-agree}(b)). Discrepancies are again where the test based on the observed score statistic rejects the null hypothesis but that based on the expected score statistic does not (Figure \ref{fig-7}, Appendix \ref{sup}).
Our simulation results for the case with $\psi_1=0.4$ (Figure \ref{fig-power04} and Table \ref{tab-agree-04} in Appendix \ref{sup}) show similar results: the power is greatly improved when the new modified rule ($T^*_O$) is used, compared to the original version of the observed score test ($T_O$). Furthermore, the power of our new test mostly exceeds the power of the other tests.

\section{Discussion}\label{sec-discuss}

After an examination of the problems with the score test statistic, 
we propose a new modified rejection rule that will work regardless of the sample size of the study. We have addressed a test of relevance also to an area of great importance in modern ecology \citep{mac06, guillera17}.

Large-sample theory, as shown for example in \citet{efron78}, supports the use of the observed information over the expected in many settings.
For the score test in our application the problem is that when $\btheta_T\ne \bM\btheta'_S$, we have seen that $E_{\btheta_T}\left(\bJ(\bM\btheta'_S)\right)$ need not be positive definite.
This is the same problem noted in the large-sample treatment in  equations (11) and (12)  of \citet{freedman07}. As in our case this leads to an ambiguous score function producing some positive and some negative eigenvalues of the observed information matrix. As a result, the observed score test statistic may be negative.

Our results for the zero-inflated Binomial model comparing two samples agree with those found for the zero-inflated Poisson model for a single parameter  in \citet{morgan07}, as well as results found in \citet{freedman07} and comments made in \citet{verbeke07}.
In \citet{morgan07}, the test was of zero inflation, while ours can be interpreted as comparing levels of zero inflation. We propose a new practical solution which makes the observed score test always useable, even when a negative test value is obtained. 

With the new test, power is restored and mostly outperforms competing tests, including the Wald test and particularly the expected score test.
In addition consistency is largely restored.

The new test presented here is easy to use, and inference is always possible. The observed score test may be obtained readily from numerical approximation and is computationally fast. Whereas the expected score test
requires derivations of expectations which in general may be difficult or not available in closed form, especially for complex models.


\begin{Appendix}\label{sec-score2}


    To examine the score statistic in more detail we adopt a standard approach.
Note that
$
    \bS_0(\btheta')=\bM^T\bS(\bM\btheta')
$
and
$
    \bJ_0(\btheta')=\bS'_0(\btheta')=\bM^T\bS'(\bM\btheta')\bM=\bM^T\bJ(\bM\btheta')\bM
$.
Define  $\bmu=E_{\btheta_T}(\bS(\bM\btheta'_S))$ and $\bSigma={\rm Cov}_{\btheta_T}\left( \bS'(\bM\btheta'_S) \right)$.
    Then
    $
    \bS(\bM\btheta'_S) \sim N(\bmu,\bSigma)
    $ approximately for large sample sizes.

Now, standard arguments yield
$
    0=\bS_0(\widehat\btheta')\approx \bS_0(\btheta'_S)+\bJ_0(\btheta'_S)\left(\widehat\btheta'-\btheta'_S \right)
$
or
$
    \left(\widehat\btheta'-\btheta'_S \right)\approx -\left(\bJ_0(\btheta'_S) \right)^{-1}\bS_0(\btheta'_S)
$.
Hence
\begin{align*}
    \bS(\bM\widehat\btheta')&\approx  \bS(\bM\btheta'_S)+\bJ(\btheta'_S)\bM(\left(\widehat\btheta'-\btheta'_S \right))\\
    &\approx \bS(\bM\btheta'_S)-\bJ(\btheta'_S)\bM \left(\bJ_0(\btheta'_S) \right)^{-1}\bS_0(\btheta'_S)\\
    &=\bS(\bM\btheta'_S)-\bJ(\btheta'_S)\bM \left(\bM^T\bJ(\bM\btheta'_S)\bM \right)^{-1}\bM^T \bS(\bM\btheta'_S)\\
    &=\left(\bI-\bJ(\btheta'_S)\bM \left(\bM^T\bJ(\bM\btheta'_S)\bM \right)^{-1}\bM^T\right)\bS(\bM\btheta'_S),
\end{align*}
where $\bI$ is the identity matrix.
Then, taking $\hat\theta'_S \approx \theta'_S$, (\ref{eq-score-2}) is approximately
\begin{align*}
    & \bS(\bM\btheta'_S)^T \left(\bI-\bM \left(\bM^T\bJ(\bM\btheta'_S)\bM \right)^{-1}\bM^T\bJ(\btheta'_S)\right)
    \left(\bJ(\bM\btheta'_S) \right)^{-1}\\
    &\times
    \left(\bI-\bJ(\btheta'_S)\bM \left(\bM^T\bJ(\bM\btheta'_S)\bM \right)^{-1}\bM^T\right)\bS(\bM\btheta'_S)\\
    &= \bS(\bM\btheta'_S)^T
    \left( \left(\bJ(\bM\btheta'_S) \right)^{-1}- \bM \left(\bM^T\bJ(\bM\btheta'_S)\bM \right)^{-1}\bM^T\right)
    \bS(\bM\btheta'_S)\\
    &= \bS(\bM\btheta'_S)^T \bSigma^{-1/2}\bSigma^{1/2}
    \left( \left(\bJ(\bM\btheta'_S) \right)^{-1}- \bM \left(\bM^T\bJ(\bM\btheta'_S)\bM \right)^{-1}\bM^T\right)
    \bSigma^{1/2}\bSigma^{-1/2}
    \bS(\bM\btheta'_S).
\end{align*}
 Define $\bP$ so that
\[
\bP^T \bSigma^{1/2}
    \left( \left(\bJ(\bM\btheta'_S) \right)^{-1}- \bM \left(\bM^T\bJ(\bM\btheta'_S)\bM \right)^{-1}\bM^T\right)
    \bSigma^{1/2} \bP=\bLambda
\]
for a diagonal matrix $\bLambda$ and $\bP^T\bP=\bP \bP^T=\bI$. Then the diagonal of $\bLambda$ is the vector of eigenvalues of
$\bSigma^{1/2}\left( \left(\bJ(\bM\btheta'_S) \right)^{-1}- \bM \left(\bM^T\bJ(\bM\btheta'_S)\bM \right)^{-1}\bM^T\right)
    \bSigma^{1/2}$ or equivalently $\left( \left(\bJ(\bM\btheta'_S) \right)^{-1}- \bM \left(\bM^T\bJ(\bM\btheta'_S)\bM \right)^{-1}\bM^T\right)
    \bSigma$.
    Thus (\ref{eq-score-2}) is approximately
    \begin{equation}
    \bS(\bM\btheta'_S)^T \bSigma^{-1/2}\bP \bLambda \bP^T \bSigma^{-1/2}
        \bS(\bM\btheta'_S).\label{eq-2}
    \end{equation}

    Write $\bS(\bM\btheta'_S)=\bmu+\Sigma^{1/2}\bZ$ where $\bZ\sim N_4(0,\bI)$. That is, $\bZ=\bSigma^{-1/2}\left(\bS(\bM\btheta'_S)-\bmu \right)$.
    Let $\bU=\bP^T \bZ$ and $\bb=\bP^T\bSigma^{-1}\bmu$. Then (\ref{eq-2}) is
    \[
    (\bb+\bU)^T \bLambda (\bb+\bU)=\sum_{j=1}^4 \lambda_j (b_j+U_j)^2.
    \]
   To examine this, we replace $ \bJ(\bM\btheta'_S)$ by $E_{\btheta_T}\left( \bJ(\bM\btheta'_S) \right)$, the expected information matrix. Note that if the null hypothesis is true, then this is $\bSigma$ and $\bmu=0$. \\
   Then $\bSigma^{1/2}\left( \bSigma^{-1}- \bM \left(\bM^T\bSigma\bM \right)^{-1}\bM^T\right)\bSigma^{1/2}$  is easily seen to be idempotent with trace 1. We now examine
    \begin{equation}
    \left( \left(E_{\btheta_T}\left( \bJ(\bM\btheta'_S) \right) \right)^{-1}- \bM \left(\bM^T( E_{\btheta_T}\left( \bJ(\bM\btheta'_S)  \right)\bM \right)^{-1}\bM^T\right)
    \bSigma,
    \end{equation}
    under our standard configuration with $\psi_1=0.8$ and take $\psi_2=\psi_1(1-R)$ for $0 \le R <1$ in steps of 0.01. The calculations in Appendix \ref{sec-moments} allow us to determine $\bSigma$ and we have previously determined  $E_{\btheta_T}\left( \bJ(\bM\btheta'_S)\right)$. We observe that, for each value of $R$, only the first eigenvalue is nonzero.
      In Figure \ref{fig-5} we plot the inverse of this nonzero eigenvalue as a function of $R$. As in our earlier examination in Section~\ref{sec-eigen} of $E_{\btheta_T}\left( \bS(\bM\btheta'_S) \right)$,
      we see that the eigenvalue becomes negative at $R \approx 0.5$. This confirms that the negative values of the score statistic are not just due to random variation.

      Clearly, if there is only one nonzero eigenvalue and this is negative then the matrix must be negative definite.
      However, the values of the score statistic were observed in our simulations to be positive and negative. To examine this, we simulate data under our standard configuration, with $\psi_1=0.8$ and $\psi_2=0.6$. This gives $\btheta'_S=(0.673,0.532,0.336)^T$. For each set of data, we compute $\bS(\bM\btheta')$
       and $\bJ(\bM\btheta')$. We take $\bSigma$ to be the empirical covariance matrix of the  $\bS(\bM\btheta'_S)$ computed for the simulated score functions and then compute the eigenvalues of $\left( \left(\bJ(\bM\btheta'_S) \right)^{-1}- \bM \left(\bM^T\bJ(\bM\btheta'_S)\bM \right)^{-1}\bM^T\right)
           \bSigma$. The inverse of the first eigenvalue is plotted in Figure \ref{fig-6}. It is apparent that the eigenvalues for the observed information matrix can be negative or positive.
           That is, random variation leads to the positive eigenvalues and hence positive values of the score statistics.

      \begin{figure}[ht]
        \centering
          \includegraphics[width=.8\textwidth]{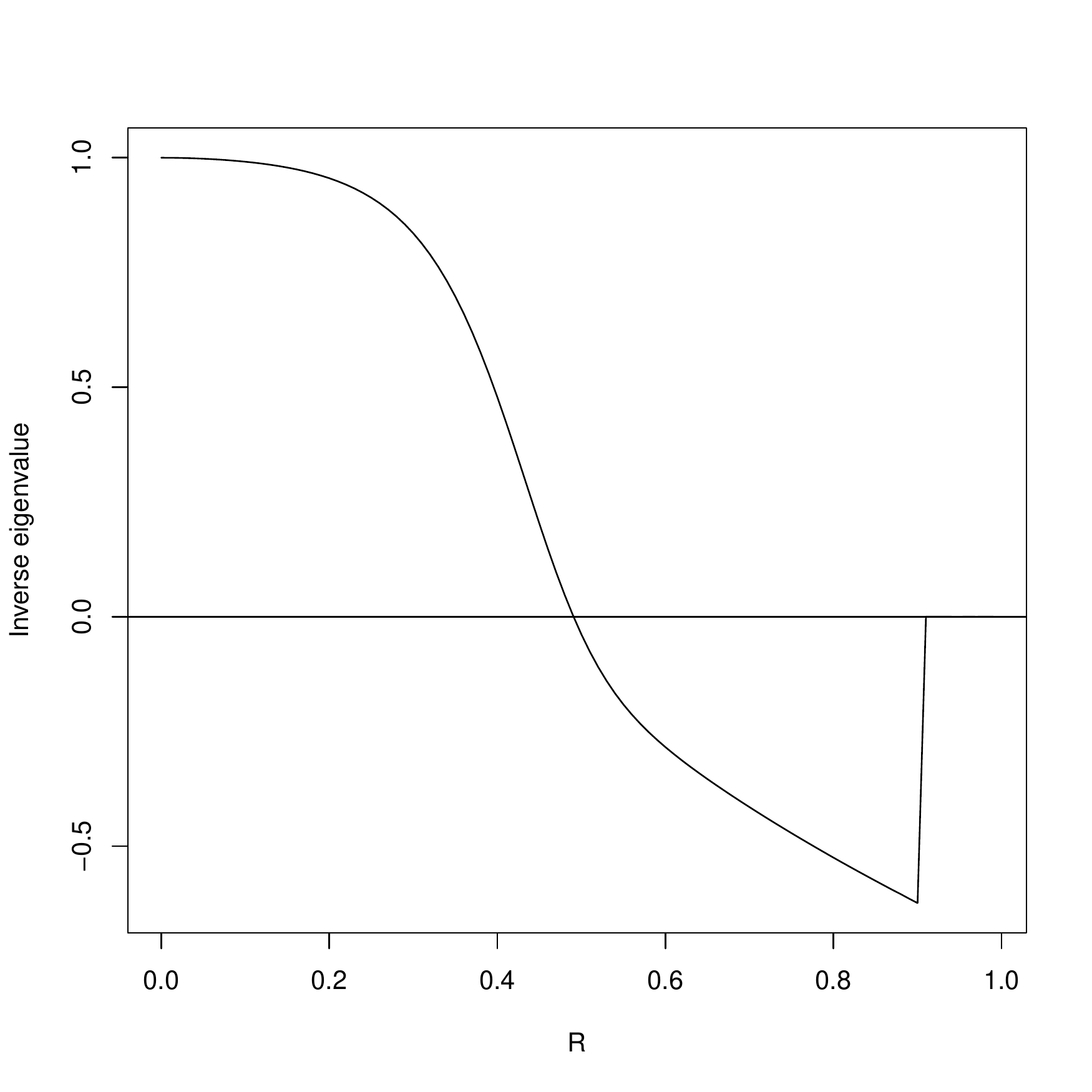}
        \caption{Inverse of the first eigenvalues of $ \left( \left(E_{\btheta_T}\left( \bJ(\bM\btheta'_S) \right) \right)^{-1}- \bM \left(\bM^T( E_{\btheta_T}\left( \bJ(\bM\btheta'_S)  \right)\bM \right)^{-1}\bM^T\right)
    \bSigma$ as a function of effect size $R$.}
        \label{fig-5}
      \end{figure}

      \begin{figure}[ht]
        \centering
          \includegraphics[width=.8\textwidth]{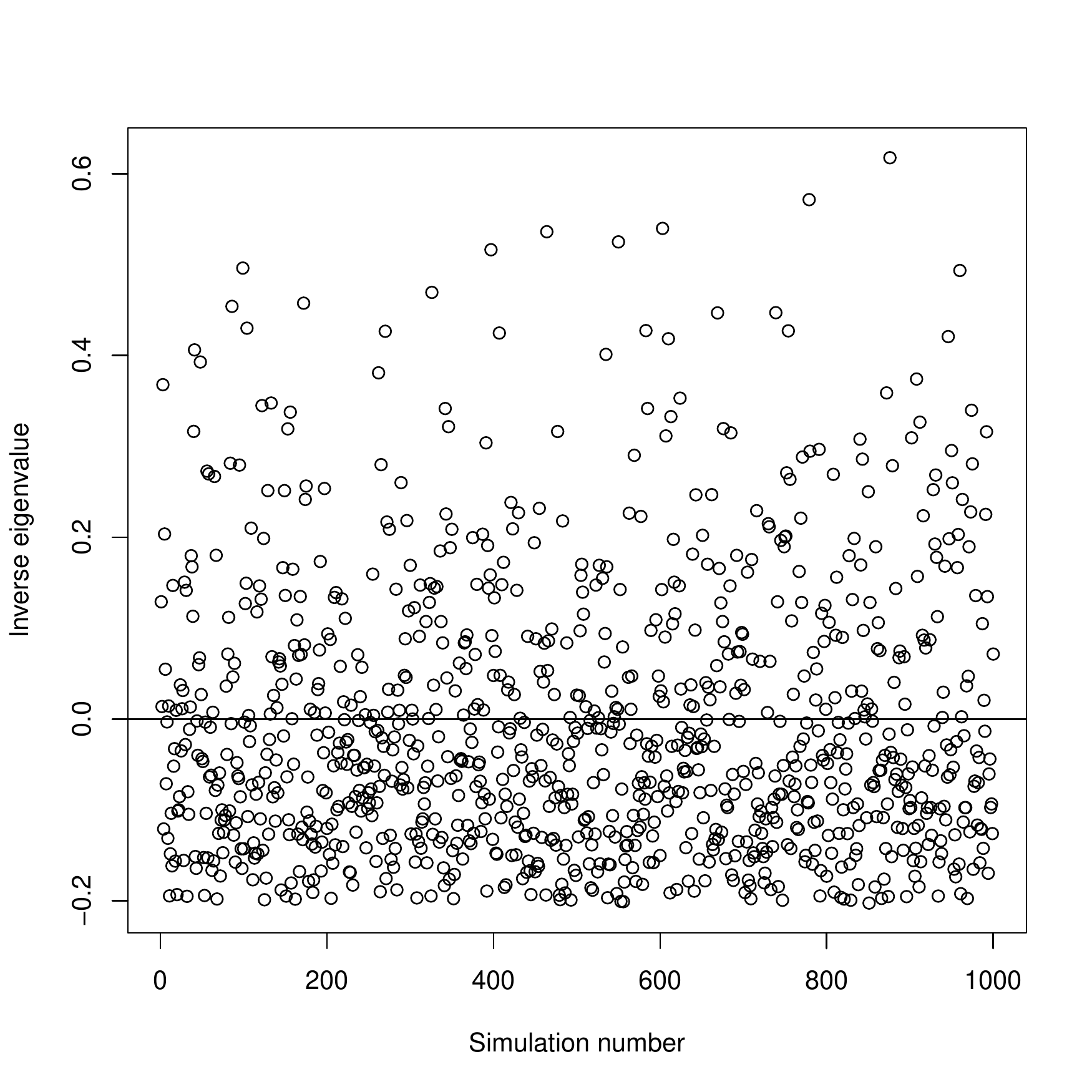}
        \caption{Inverse of the first eigenvalue of $\left( \left(\bJ(\bM\btheta'_S) \right)^{-1}- \bM \left(\bM^T\bJ(\bM\btheta'_S)\bM \right)^{-1}\bM^T\right)
           \bSigma$ when $R=0.6$, for 1000 simulations.}
        \label{fig-6}
      \end{figure}

\end{Appendix}

\begin{Appendix}\label{sec-moments}

To determine the eigenvalues we need to compute $\bSigma$.
Recall that for a single region,
\begin{align*}
    S_1&=\frac{s_d-\psi \theta N}{\psi(1-\psi\theta)},\\
    S_2&= \frac{d-s_d K p}{p(1-p)}-\frac{(N-s_d)\psi K (1-\theta)}{(1-\psi\theta)(1-p)},
\end{align*}
$s_d \sim {\rm Bin}(N,\psi_T\theta_T)$ and given $s_d$, $d$ is the sum of $s_d$ independent positive binomial random variables.
Then
\begin{align*}
E(s_d)&=N\psi_T\theta_T,\\
{\rm Var}(s_d)&=N\psi_T\theta_T(1- \psi_T\theta_T),\\
E(d\vert s_d)&= s_d \frac{K p_T}{\theta_T},\\
{\rm Var}(d\vert s_d)&=s_d \left( \frac{K^2p_T^2-K p_T^2+K p}{\theta_T}-\frac{K^2 p_T^2}{\theta_T^2}\right), \\
 E(d)&=N\psi_TK p_T,\\
 {\rm Var}(d)&=E({\rm Var}(d\vert s_d))+{\rm Var}(E(d\vert s_d)),\\
 &=N\psi_T\theta_T\left( \frac{K^2p_T^2-K p_T^2+K p}{\theta_T}-\frac{K^2 p_T^2}{\theta_T^2}\right)+
 N\psi_T\theta_T(1- \psi_T\theta_T) \left(\frac{K p_T}{\theta_T} \right)^2,\\
 {\rm Cov}(d,s_d)&=E(d s_d)-E(s_d)E(d),\\
 &=E(s_d^2 K p_T/\theta_T)- N^2\psi_T^2\theta_T K p_T,\\
 &= \left(N\psi_T\theta_T(1- \psi_T\theta_T)+N^2\psi_T^2\theta_T^2 \right) K p_T/\theta_T- N^2\psi_T^2\theta_T K p_T\\
 &=N\psi_T(1-\psi_T\theta_T)Kp_T.
\end{align*}

Thus
\begin{align*}
    \bmu_1&= \frac{N(\psi_T\theta_T-\psi\theta)}{\psi(1-\psi\theta)}.\\
    \bmu_2&=\frac{N\psi_T K\left(p_T-p\theta_T\right)}{p(1-p)}-\frac{N(1-\psi_T\theta_T)\psi K(1-\theta)}{(1-\psi\theta)(1-p)}.
\end{align*}

\begin{align*}
 \bSigma_{11}={\rm Var}(\bS_1)&= \frac{N\psi_T\theta_T(1- \psi_T\theta_T)} {\psi^2(1-\psi\theta)^2},\\
  \bSigma_{22}={\rm Var}(\bS_2)&= E({\rm Var}(S_2\vert s_d))+{\rm Var}(E(S_2\vert s_d))\\
 &=\frac{ E\left({\rm Var}(d\vert s_d)\right)}{p^2(1-p)^2}+{\rm Var}(s_d)\left(\frac{Kp_T-Kp\theta_T}{\theta_Tp(1-p)}+\frac{\psi K(1-\theta)}{(1-\psi\theta)(1-p)} \right)^2,
\end{align*} 
\begin{align*} 
  E(\bS_1\bS_2)&=E\left\{(s_d-\psi\theta N)s_d\right\}\frac{Kp_T/\theta_T-Kp}{\psi(1-\psi\theta)p(1-p)} \\
  &-E\left\{(s_d-\psi\theta N)(N-s_d) \right\}\frac{\psi K(1-\theta)}{\psi(1-\psi\theta)^2(1-p)}\\
  &=\left\{E(s_d^2)-\psi\theta N E(s_d)\right\}\frac{Kp_T/\theta_T-Kp}{\psi(1-\psi\theta)p(1-p)}\\
  &- \left\{N(E(s_d)-\psi\theta N)-E(s_d^2)+\psi\theta N E(s_d) \right\}\frac{\psi K(1-\theta)}{\psi(1-\psi\theta)^2(1-p)},\\
  \bSigma_{12}&= E(\bS_1\bS_2)-\bmu_1\bmu_2.
\end{align*}
\end{Appendix}

\newpage

\begin{Appendix}\label{sup}

\begin{figure}[h!]
\begin{center}
\includegraphics[width=.8\textwidth]{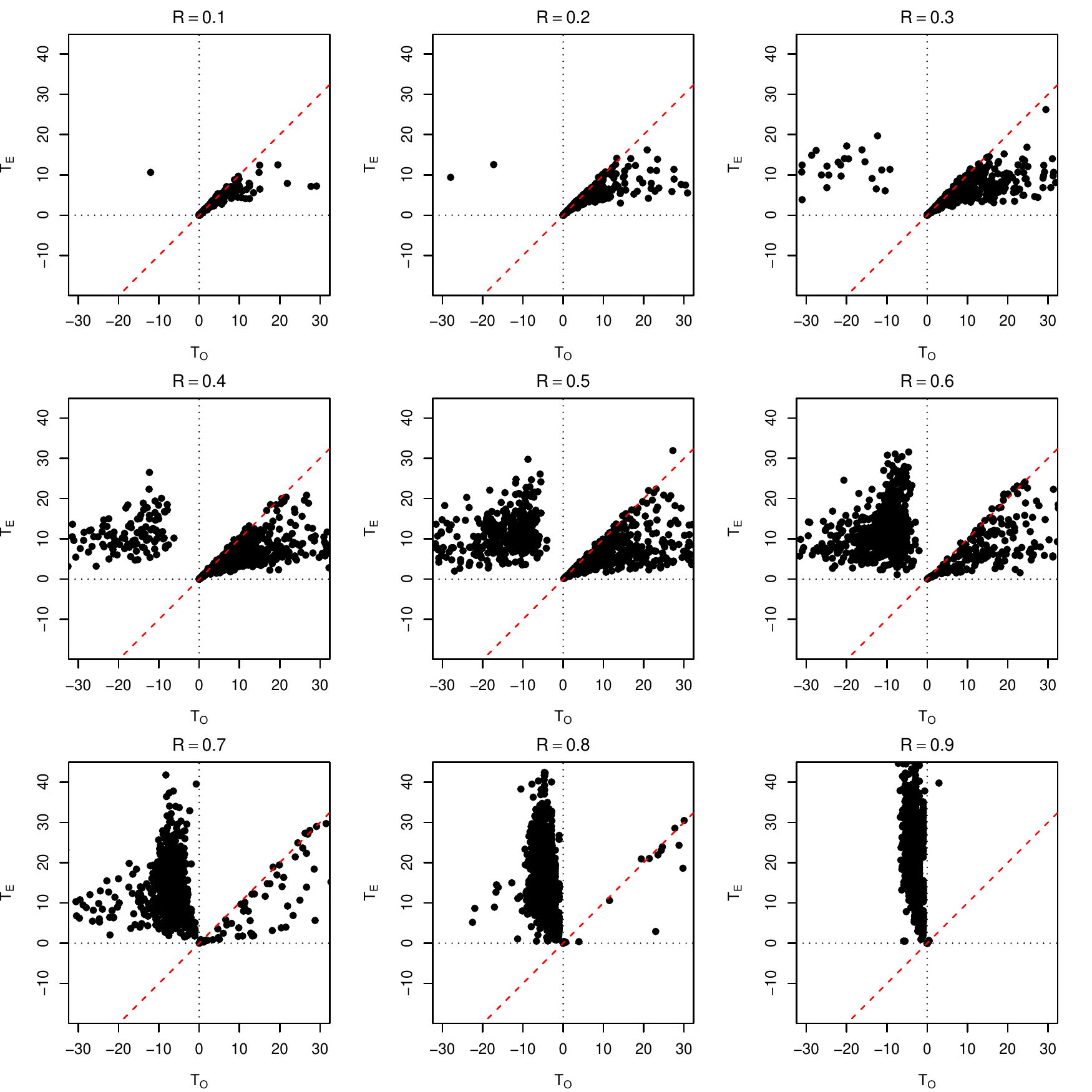}
\end{center}
\caption{Observed ($T_O$) versus expected ($T_E$) score test statistic, for $\psi_1 = 0.8$ under the standard configuration. For clarity, 1000 of the 50000 simulations are displayed here and axes limits are set to $\pm 30$.}
\label{fig-7}
\end{figure}

\begin{figure}[h!]
\begin{center}
\includegraphics[width=.8\textwidth]{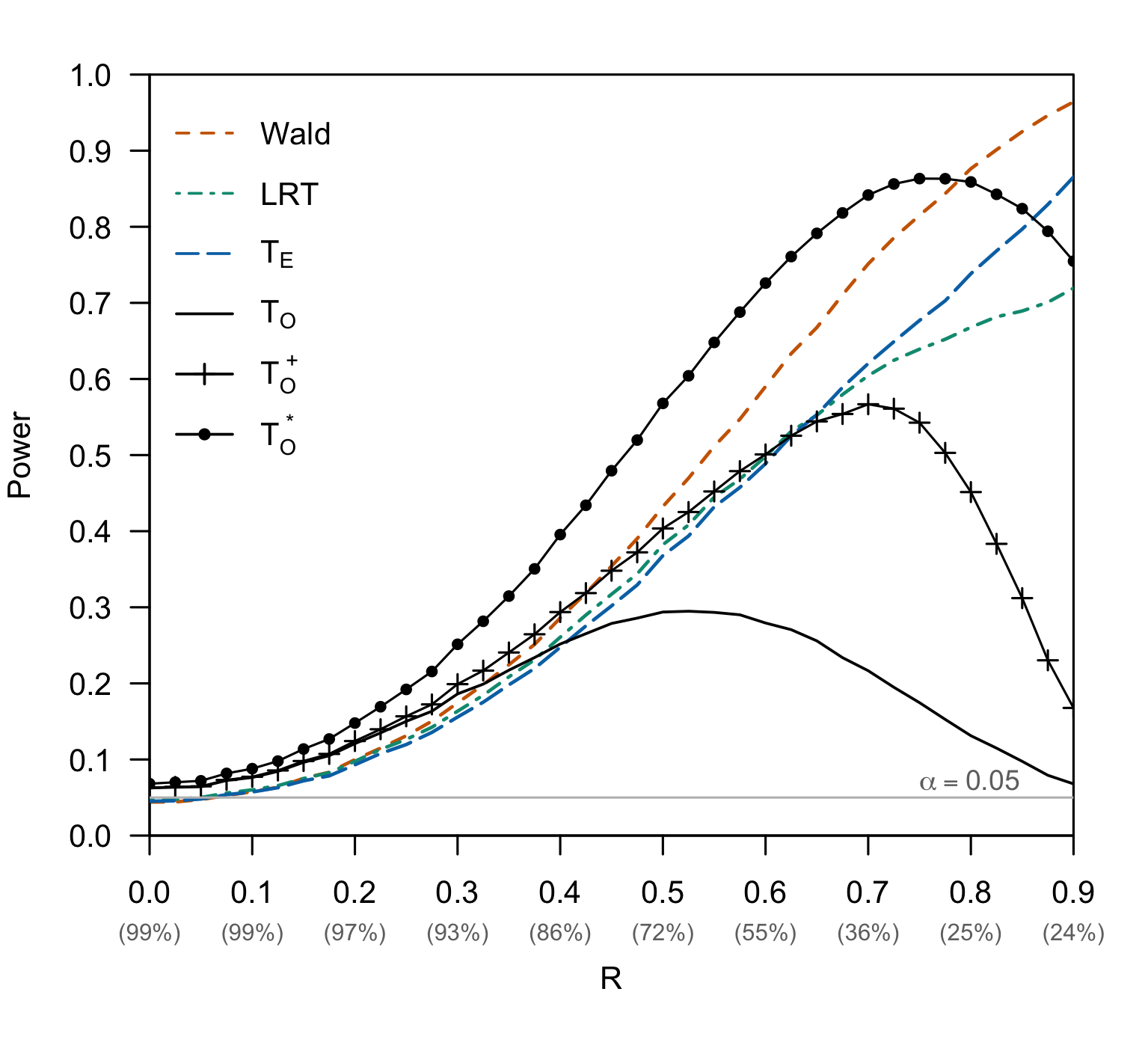}
\end{center}
\caption{Power plot for scenario $\psi_1=0.4$ under the standard configuration, for different effect sizes $R$ and computed via simulation (50000 datasets per $R$ point). Each curve represents one of the tests considered (LRT = Likelihood-ratio test; Wald = Wald test; $T_E$ = Score test using expected information; $T_O$ = Score test using the observed information; $T_O^*$ = New modified observed score test; $T_O^+$ = Observed score test when only positive test statistics are considered). Results are displayed for simulations where numerical optimization did not fail. Percentage of simulations that returned a positive observed score test ($T_O^+$) are listed in parenthesis under the $x$-axis.}
\label{fig-power04}
\end{figure}

\clearpage

\begin{table}[ht]
    \caption{Agreement (Agr.) for accepting/rejecting $H_0$ between tests based on the observed  and expected score statistics, for data simulated under the standard configuration with $\psi_1=0.4$. In a) only datasets where the observed score statistic is positive ($T_O^+$) are considered. In b) both positive and negative score statistics are considered under our new modified rejection rule ($T_O^*$). In each case, $n$ is the number of simulated tests considered (out of the 50000 simulations). As effect size $R$ increases, there are fewer simulations with a positive observed score test statistic, hence the corresponding decrease in $n$.}
\vspace{2mm}
\centering
\begin{small}
\begin{tabular}{|l|l|rrrrrrrrrr|}
  \hline
& $R$ & 0 & 0.1 & 0.2 & 0.3 & 0.4 & 0.5 & 0.6 & 0.7 & 0.8 & 0.9 \\
  \hline
a) $T_O^+$ & Agr.\ & 0.97 & 0.97 & 0.95 & 0.92 & 0.88 & 0.85 & 0.82 & 0.84 & 0.92 & 0.99 \\ \hline
 Positive &  $n$ & 49721 & 49423 & 48648 & 46729 & 42786 & 36203 & 27412 & 18123 & 12380 & 11894 \\
 \hline \hline
  b) $T_O^*$ & Agr.\ & 0.97 & 0.97 & 0.94 & 0.90 & 0.85 & 0.80 & 0.75 & 0.75 & 0.80 & 0.91 \\ \hline
  Modified &  $n$ & 50000 & 50000 & 49999 & 49997 & 49993 & 49976 & 49932 & 49593 & 48137 & 40369 \\
   \hline
\end{tabular}
\end{small}
\label{tab-agree-04}
\end{table}

\end{Appendix}

\newpage
\bibliography{HT-bib-V6-NK}
\bibliographystyle{apalike}

\end{document}